\documentclass[12pt,twoside]{article}
\usepackage[english]{babel}
\usepackage[latin1]{inputenc}
\usepackage{amsmath}
\usepackage{mathtools} 
\usepackage{cancel}
\usepackage{tikz}
\usetikzlibrary{calc}

\usepackage{graphicx}
\usepackage{times,amssymb,amscd}
\newtheorem{thm}{Theorem}[section]
\newtheorem{cor}[thm]{Corollary}

\newtheorem{defn}[thm]{Definition}
\newtheorem{rem}[thm]{Remark}
\numberwithin{equation}{section}

\newcommand{\bH}{\mathbf{H}}

\newcommand{\bL}{\mathbf{L}}

\newcommand{\bR}{\mathbf{R}}
\newcommand{\bS}{\mathbf{S}}

\newcommand{\bx}{\mathbf{x}}
\newcommand{\by}{\mathbf{y}}

\newcommand{\bu}{\mathbf{u}}
\newcommand{\bv}{\mathbf{v}}

\newcommand{\rS}{\mathrm{S}}
\newcommand{\rC}{\mathrm{C}}
\newcommand{\rT}{\mathrm{T}}

\begin{document}
\pagestyle{myheadings}
\markboth{\centerline{G\'eza Csima }}
{Generalized Apollonius Circles As Equioptic Curves $\dots$}
\title
{Generalized Apollonius Circles As Equioptic Curves Of Circles In Constants Curvature Geometries
}


\author{G\'eza Csima \\
\normalsize Department of Algebra and Geometry, Institute of Mathematics,\\
\normalsize Budapest University of Technology and Economics, \\
\normalsize M\H uegyetem rkp. 3.,\\
\vspace{0.5cm}\normalsize H-1111 Budapest, Hungary \\
\normalsize Department of Basic Technical Studies,\\
\normalsize University of Debrecen,\\
\normalsize \'Otemet\H o str. 2--4., \\
\normalsize H-4028 Debrecen, Hungary\\
\normalsize csimageza@gmail.com\\
\date{\normalsize{\today}}}

\maketitle
\begin{abstract}
We extend the old definition of the Apollonius circle in such a way that it results in the same curve in Euclidean geometry but will be more convenient in hyperbolic and spherical geometries. We show that there exists an Apollonius circle of the centers of two circles that coincides with their equioptic curves, as in Euclidean geometry.
\end{abstract}

\section{Introduction} \label{section1}

Apollonius of Perga asked for the first time: \textit{"What is the locus of points for which the ratio of the distances to two given points is constant k?"} The solution can be considered as an alternative definition of a circle for $k\neq1$. Many applications are connected to this topic, Apollonian packing, e.g. It is also interesting that the isoptic curves of a segment are the orthogonal trajectories of the Apollonian hyperbolic circle pencil, defined by the endpoints of the segment. The bipolar coordinate system is based on this property. Now, we consider an interesting problem from the work of Odehnal \cite{Odehnal}, where the result is the Apollonian circle.
\begin{defn}
    Given two plane curves $c_1$ and $c_2$ we call the set of points from which $c_1$ and $c_2$ are seen under equal angle the equioptic curve.
    \label{def:equiopt}
\end{defn}
In Section 3.1, it is proved that if $c_1$ and $c_2$ are circles, then their equioptic curve is either a line (congruent case), empty (concentric case), or a circle (general case) regarding the real (not ideal) points. It is easy to see that the diameter of this circle is defined by the segment connecting the centers of similitude. 
\begin{cor}
    The equioptic curve of two circles with radius $r_1$ and $r_2$ is an Apollonius circle of the centers with $k=\dfrac{r_1}{r_2}.$ 
    \label{cor:App_equiopt}
\end{cor}

\textit{The goal of this study is to prove a statement in hyperbolic and elliptic geometries similar to Corollary \ref{cor:App_equiopt}.}

\section{Projective model}
The following model (see \cite{MSz}) is capable of uniting the 8 Thurston geometries (\cite{T}). Let $G$ be one of the constant curvature plane geometries. 
We use homogeneous coordinates $\mathbf{x} = (x^1 : x^2 : x^3)$ and $\mathbf{u} = (u_1 : u_2 : u_3)$
in order to represent points $X$ as well as straight lines $u$ in projective space $\mathbb{P}^3$.
Sometimes we write $X = \mathbf{x}\mathbb{R}$ in order to express that the point $X$ is determined by
a certain vector $\mathbf{x}\in\mathbb{R}^3$ and its non-trivial scalar multiples, similarly $u = \mathbf{u}\mathbb{R}$ for
lines. A point $X =\mathbf{x}\mathbb{R}$ and a straight line $u=\mathbf{u}\mathbb{R}$ are incident if $\mathbf{x}\cdot\mathbf{u}^T = 0$
with $\cdot$ denoting the usual matrix multiplication.

Constant curvature plane geometries can be represented in projective space $\mathbb{P}^3$ using the bilinear form
\begin{equation}\label{222.1}
\langle \mathbf{x}, \mathbf{y} \rangle = x^1y^1+x^2y^2+\epsilon x^3 y^3, 
\end{equation}
where $\epsilon=0,-1,+1$, respectively, Euclidean, hyperbolic and elliptic geometries. Now we consider hyperbolic and elliptic planes simultaneously, so hereafter $\epsilon=\pm1,$ $+1$ for elliptic an $-1$ for hyperbolic geometry.

Hyperbolic and elliptic distances and angles can be computed with the help of the bilinear form \eqref{222.1}, see \cite{MSz}. 
The distance $d(X,Y)$ of two proper points $X = \bx\mathbb{R}$ and $Y =\by\mathbb{R}$ is
given by
\begin{equation}\label{222.2}
\mathrm{C(d(X,Y))}=\frac{\epsilon\langle \mathbf{x},\mathbf{y} \rangle }{\sqrt{\langle \mathbf{x},\mathbf{x} \rangle
\langle \mathbf{y},\mathbf{y} \rangle }},
\end{equation}
where $\mathrm{C(*)}$ is the cosine function in elliptic geometry and the hyperbolic cosine function in hyperbolic geometry. We also introduce the $\mathrm{S(*)}$ and $\rT(*)$ functions, which are the sine and tangent functions in elliptic geometry and the hyperbolic sine and tangent functions in hyperbolic geometry. 


Finally, we find the angle $\alpha(u,v)$ enclosed by two proper
straight lines $u=\bu\mathbb{R}$ and $v=\bv\mathbb{R}$, respectively, with
\begin{equation}
\cos{\alpha}=\frac{\epsilon\langle \bu,\bv \rangle}{\sqrt{
\langle \bu,\bu\rangle \langle \bv,\bv \rangle}}.
\label{222.3}
\end{equation}

\section{Apollonian circles}
There are surprisingly few results regarding Apollonian circles in hyperbolic geometry but there are plenty regarding Apollonian packings and Apollonius problem. In \cite{Ion} Ione\c{s}cu determines the hyperbolic equioptic curves of consecutive collinear segments that is an equivalent approach to Apollonian circles in the Euclidean geometry. The results are generally quartic curves in polar form in the half-plane model. 

\subsection{First approach}
In the elliptic/hyperbolic plane, the most straightforward method to define the Apollonius circle is by replacing Euclidean distance
with elliptic/hyperbolic distance, described in $(\ref{222.2}).$

Let
\[
A=(a:0:1), \qquad B=(b:0:1),
\]
be two points on the $x$-axis and let
\[
P=(x:y:1)
\]
be an arbitrary point in the disk. Substituting into the defining Apollonius ratio yields
\[
\rC^{-1}\!\left(
\frac{1+\epsilon\,ax}
{\sqrt{(1+\epsilon\,a^2)(1+\epsilon(x^2+y^2))}}
\right)
=
k \cdot
\rC^{-1}\!\left(
\frac{1+\epsilon\,bx}
{\sqrt{(1+\epsilon\,b^2)(1+\epsilon(x^2+y^2))}}
\right),
\]
that is a nonlinear (Chebyshev-like) implicit equation for the
Apollonius locus. To avoid this, we generalize the definition of the Apollonian circle.

\subsection{Generalized Apollonian curve}

We will deal with the hyperbolic and elliptic cases together. It is known that the circumference of a hyperbolic/elliptic circle is $C_r=2\pi \cdot\mathrm{S}(r)$ where $r$ is the radius and $\mathrm{S(*)}$ is the sine function in elliptic geometry and hyperbolic sine function in hyperbolic geometry.
\begin{defn}
    Let us be given two points $A$ and $B$ and $k\in\mathbb{R}^+.$ The locus of points $P$ from which the ratio of the circumference of the circles around $P$ to $A$ and $B$ is $k$ are called generalized Apollonian curve. 
    \label{def:Ap_gen}
\end{defn}
It is easy to see that this paraphrase will not change the obtained curve in Euclidean geometry but in spherical and hyperbolic geometry:
\begin{equation}
    \dfrac{C_{|\overline{AP}|}}{C_{|\overline{BP}|}}=\dfrac{\bcancel{2\pi}\cdot\mathrm{S}(d(A,P))}{\bcancel{2\pi}\cdot\mathrm{S}(d(B,P))} = k.
\end{equation}

For the sake of nice results and easy computations, we transform the problem in the coordinate system to the most convenient position by isometries. Without loss of generality, we can assume that $A(a:0:1)$ and $B(b:0:1)$ points are situated on the $x$ axis such that the origin $O(0:0:1)$ is between them and $a>0$. Then $b<0$ follows. Furthermore, we can also assume that $O$ is part of the generalized Apollonian curve: $\mathrm{S}(d(A,O))=k\cdot\mathrm{S}(d(B,O))$ 
It is easy to compute that if $d(A,B)=d,$ then 
\begin{equation}
    a=\dfrac{k\,\rS(d)}{1+k\,\rC(d)}\ \mathrm{and}\ b=-\dfrac{\rS(d)}{\rC(d)+k}.
\end{equation} 
We note here that in the elliptic case due to the nature of the geometry we assume that $d<\dfrac{\pi}{2}.$ \\
With the well known $\rC^2(t)+\epsilon\,\rS(t)^2=1$ identity and (\ref{222.2}), we know that

\footnotesize
\begin{equation}
\begin{gathered}
    \rS^2(d(A,P))=\epsilon\,(1-\rC^2(d(A,P)))=\epsilon\left(1-\dfrac{(1+\epsilon ax)^2}{(1+\epsilon\,(x^2+y^2))(1+\epsilon\,a^2)}\right)= \\
    \dfrac{\bcancel{1}+\epsilon\,a^2+\epsilon\,x^2+\bcancel{a^2x^2}+y^2(a^2+\epsilon)-\bcancel{1}-2\epsilon\,ax+\bcancel{a^2x^2}}{(1+\epsilon\,(x^2+y^2))(a^2+\epsilon)}=\dfrac{y^2(a^2+\epsilon)+\epsilon\,(x-a)^2}{(1+\epsilon\,(x^2+y^2))(a^2+\epsilon)}
\end{gathered}
\end{equation}
\normalsize
Similarly, we can obtain $\rS^2(d(B,P))$ by replacing $a$ with $b.$ Then after canceling $(1+\epsilon\,(x^2+y^2))$ we obtain 
\begin{equation}
    k^2=\dfrac{\mathrm{S}^2(d(A,P))}{\mathrm{S}^2(d(B,P))} = \dfrac{(y^2(a^2+\epsilon)+\epsilon\,(x-a)^2)(b^2+\epsilon)}{(y^2(b^2+\epsilon)+\epsilon\,(x-b)^2)(a^2+\epsilon)}=\dfrac{y^2+\dfrac{(x-a)^2}{1+\epsilon\,a^2}}{y^2+\dfrac{(x-b)^2}{1+\epsilon\,b^2}}
\end{equation}
Ordering this equation, we obtain that
\begin{equation}
    \dfrac{(x-a)^2}{1+\epsilon\,a^2}-k^2\cdot\dfrac{(x-b)^2}{1+\epsilon\,b^2}+(1-k^2)y^2=0,
\end{equation}
that is a quadratic curve. It is easy to see that $\dfrac{1}{1+\epsilon\,a^2}-\dfrac{k^2}{1+\epsilon\,b^2}$ is the coefficient of $x^2$ but according to our assumption, the origin must satisfy this equation. Therefore, replacing $x$ and $y$ with 0, we get that $\dfrac{a^2}{1+\epsilon\,a^2}-\dfrac{k^2b^2}{1+\epsilon\,b^2}$ must be equal to $0.$ Then we can see that
\begin{equation}
\begin{gathered}
    \dfrac{a^2}{1+\epsilon\,a^2}=\dfrac{k^2b^2}{1+\epsilon\,b^2}\Leftrightarrow \dfrac{\epsilon\,a^2}{1+\epsilon\,a^2}=\dfrac{\epsilon\,k^2b^2}{1+\epsilon\,b^2}\Leftrightarrow
    \dfrac{1+\epsilon\,a^2-1}{1+\epsilon\,a^2}=\dfrac{k^2+\epsilon\,k^2b^2-k^2}{1+\epsilon\,b^2}\\\
    \Leftrightarrow
    1-\dfrac{1}{1+\epsilon\,a^2}=k^2-\dfrac{k^2}{1+\epsilon\,b^2}\Leftrightarrow\dfrac{1}{1+\epsilon\,a^2}-\dfrac{k^2}{1+\epsilon\,b^2}=1-k^2
\end{gathered}
\end{equation}
We can summarize our results in the following theorem.
\begin{thm}
Let us be given $k,d\in\mathbb{R}^+$ and $A(a:0:1),$ $B(b:0:1)$ points in the elliptic or hyperbolic plane such that $  a=\dfrac{k\,\rS(d)}{1+k\,\rC(d)}\ \mathrm{and}\ b=-\dfrac{\rS(d)}{\rC(d)+k},$ where $k\neq1$ and in the elliptic case $d<\dfrac{\pi}{2}.$ Then the generalized Apollonian curves of $A$ and $B$ with ratio $k$ has an equation:
\begin{equation}
    x^2-2\dfrac{a(1+\epsilon\,b^2)-k^2b(1+\epsilon\,a^2)}{(1-k^2)(1+\epsilon\,a^2)(1+\epsilon\,b^2)}x+y^2=0
    \label{eq:app}
\end{equation}
\end{thm}
\begin{rem}
\begin{enumerate}
    \item The (\ref{eq:app}) equation can be given with $d$ as well
\begin{equation}x^2-2\dfrac{k\,\rS(d)
}{1-k^2}\,x+y^2=0
\end{equation}
\item Although the resulted curves are quadratic curves and they look like Euclidean circles in the model, they are not elliptic or hyperbolic circles but other type of conics. In the hyperbolic case, the result is either an ellipse(0) or a paracycle(1) or a semi-hyperbola(2) depending on the number of common points with the absolute indicated in the parentheses. See \cite{Cs-Sz16-1} for further types of hyperbolic conic sections. In the elliptic case, all conic sections are isometric images of ellipses, only the number of common points with the ideal line distinguish them.
\item For fixed $d,$ the resulted circle pencil is parabolic such that their common point is the origin $O(0:0:1)$ and their common tangent is the $x=0$ line. 
\item If $k=1$ then $a=-b$ and the result is the $x=0$ line.
\end{enumerate}
\end{rem}

\begin{figure}[h!]
    \centering
    \includegraphics[width=0.45\linewidth]{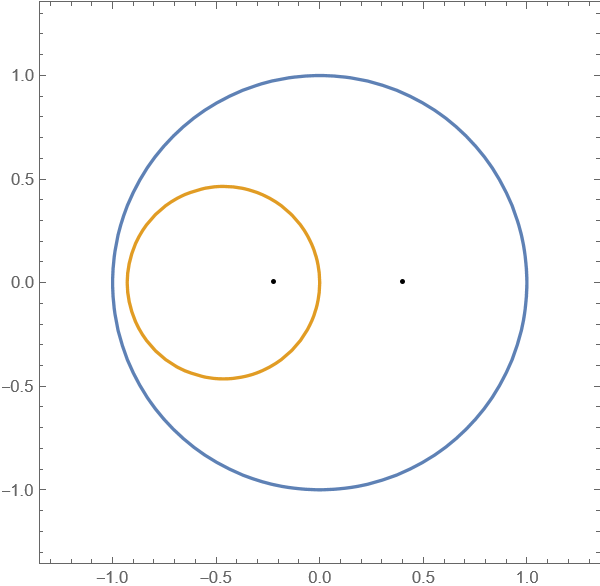} \includegraphics[width=0.45\linewidth]{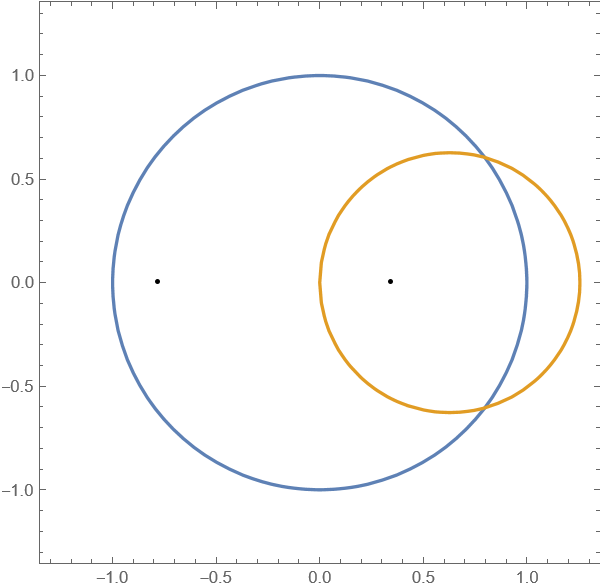}
    \caption{Generalized Apollonian curves on hyperbolic plane with $d=0.65,$ $k=2$ on the left and $d=1.4,$ $k=0.3$ on the right.}
    \label{fig:1}
\end{figure}

\begin{figure}
    \centering
    \includegraphics[width=0.45\linewidth]{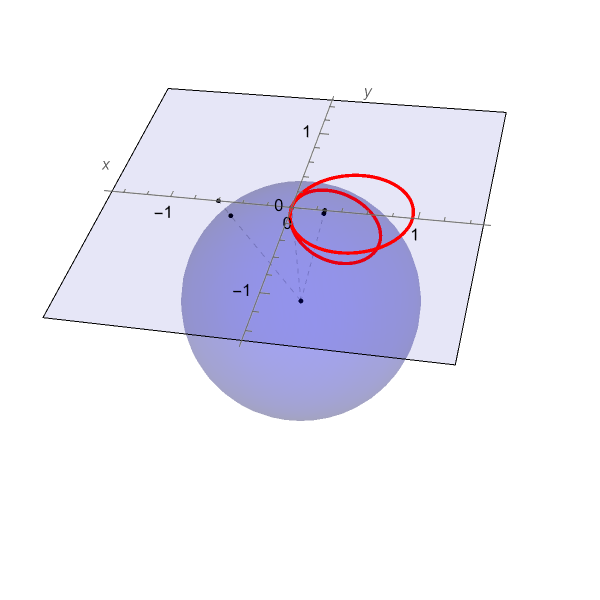} \includegraphics[width=0.45\linewidth]{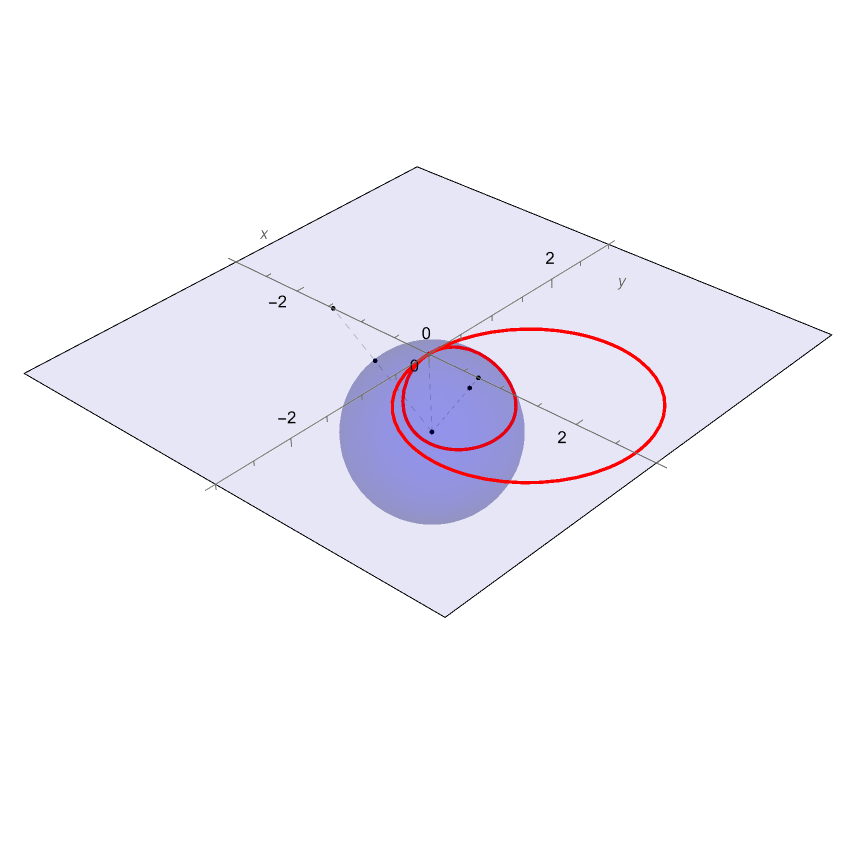}
    \caption{Generalized Apollonian curves on elliptic plane with $d=0.8,$ $k=0.5$ on the left and $d=\pi/2,$ $k=0.7$ on the right.}
    \label{fig:2}
\end{figure}

On Figure \ref{fig:1} and Figure \ref{fig:2} some generalized Apollonian curves can be seen on elliptic and hyperbolic plane.

\section{Equioptic curves}
In this section, we would like to compute the equioptic curves of elliptic/hyperbolic circles. 



\begin{figure}[h]
    \centering
    \begin{tikzpicture}[scale=1]

    \def\r{2}

    \coordinate (O) at (0,0);

    \coordinate (P) at (5,1);

    \draw[thick] (O) circle (\r);

    \fill (P) circle (2pt);
    \node[right] at (P) {$P$};

    \fill (O) circle (2pt);
    \node[below left] at (O) {$C$};

    \pgfmathsetmacro{\d}{sqrt(5^2+1^2)}

    \pgfmathsetmacro{\theta}{atan2(1,5)}
    \pgfmathsetmacro{\phi}{acos(\r/\d)}

    \coordinate (T1) at ({\r*cos(\theta+\phi)},
                          {\r*sin(\theta+\phi)});

    \coordinate (T2) at ({\r*cos(\theta-\phi)},
                          {\r*sin(\theta-\phi)});

    \draw[dashed] (O) -- (T1);
    \draw[dashed] (O) -- (T2);

    \draw[thick] (P) -- (T1);
    \draw[thick] (P) -- (T2);

    \fill (T1) circle (2pt);
    \fill (T2) circle (2pt);

    \node[above right] at (T1) {$T_1$};
    \node[below right] at (T2) {$T_2$};
    \node[above left] at (P) {$\frac{\alpha}{2}$};

    \node[left] at ($(O)!0.5!(T1)$) {$r$};

    \draw ($(T1)!0.12!(O)$);

    \draw ($(T2)!0.12!(O)$);
   \draw (P) -- (O);

\draw (T1) -- (T2);
\node[above right] at ($(O)!0.5!(P)$) {$d$};
\draw
  ($(P)+(168.69:0.6)$)
  arc[start angle=168.69,end angle=191.31,radius=0.6];
\end{tikzpicture}
    \caption{Tangent from outer point}
    \label{fig:tan}
\end{figure}
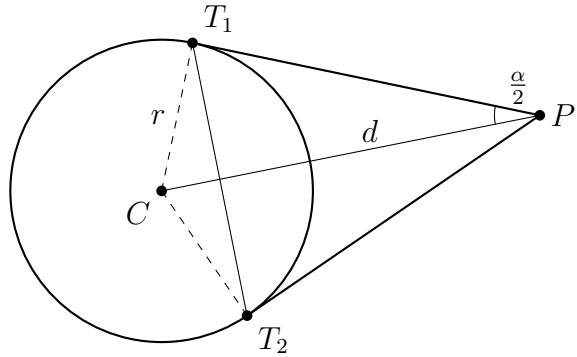

 Let us consider two circles $c_1$ and $c_2$ on the elliptic or hyperbolic plane with centers $C_{1,2}$ and radii $r_{1,2}.$
To obtain their equioptic curve, we must determine, from an arbitrary point $P,$, under what angle $\alpha_{1,2}$ the circles can be seen. We could obtain that from the results of \cite{Cs-Sz14} as a special case of ellipse, where $f=0\Rightarrow a=b,$ but we provide here a more suitable and model-free proof. 

Look at $PCT_1\triangle$ on Figure \ref{fig:tan}. We can use the sine rule for $CT_1$ and $CP$ sides with their opposing angles:
\begin{equation}
    \dfrac{\rS(d)}{\sin CT_1P\angle}=\dfrac{\rS(r)}{\sin T_1CP\angle}\Rightarrow\dfrac{\rS(d)}{1}=\dfrac{\rS(r)}{\sin\frac{\alpha}{2}}\Rightarrow\sin\frac{\alpha}{2}=\dfrac{\rS(r)}{\rS(d)}.
\end{equation}
We can apply this to both circles and if $\alpha_1=\alpha_2$, then we obtain the equioptic:
\begin{equation}
    \dfrac{\rS(r_1)}{\rS(d_1)}=\sin\frac{\alpha_1}{2}=\sin\frac{\alpha_2}{2}=\dfrac{\rS(r_2)}{\rS(d_2)}\Rightarrow\dfrac{\rS(d_1)}{\rS(d_2)}=\dfrac{\rS(r_1)}{\rS(r_2)}\Rightarrow\dfrac{\rS(d(P,C_1))}{\rS(d(P,C_2))}=k.
\end{equation}

We note here that $\alpha_i<\pi\Rightarrow\frac{\alpha}{2}<\frac{\pi}{2},$ where sine is unique. Now we can phrase the following theorem:
\begin{thm}
  Let us be given two circles $c_1$ and $c_2$ on the elliptic or hyperbolic plane with centers $C_{1,2}$ and radii $r_{1,2}.$ Then their equioptic curve is their generalized Apollonian curve with $k=\dfrac{\rS(r_1)}{\rS(r_2)}.$
\end{thm}
\section{Conclusion and future plans}
The most important part of this study is the alternative definition of Apollonian curves. The modified definition is invariant in the Euclidean plane, but in other constant curvature geometries, the results are quadratic curves and not nonlinear equations. Furthermore, we have seen that the nice property (see Corollary \ref{cor:App_equiopt}), that connects equioptic curves of circles with Apollonian circles in Euclidean geometry is also valid with this generalized Apollonian curve in elliptic/hyperbolic geometries.

Our next plan is to extend these results into 3 dimensional non-Euclidean so called Thurston geometries (see \cite{T}). For $\mathbf{H}^3$ and $\mathbf{S}^3,$ the analogous extension of Definition \ref{def:Ap_gen} to surface area of spheres seems straightforward. The question for product geometries ($\bS^2\!\times\!\bR$ and $\bH^2\!\times\!\bR$) and for twisted geometries $\widetilde{\bS\bL_2\bR},\ \mathbf{Nil}\ \mathrm{and}\ \mathbf{Sol}$ seems more challenging. In these topic there are already some results, see \cite{Cs-Sz25,Sz}.  
\section*{Data Availability Statement}
Data sharing not applicable to this article as no databases were generated or analyzed during the current study.

\section*{Conflict Of Interest Statement}
The author have no conflict of interest to declare that are relevant to the content of this study.


\begin{thebibliography}{12}
%

\bibitem{Cs-Sz14}
{{Csima,~G.~--~Szirmai,~J.:}
Isoptic curves of conic sections in constant curvature geometries.
\emph{Mathematical Communications}
{\bf 19} (2014) 277--290.}
%
\bibitem{Cs-Sz16-1}
{Csima,~G.~--~Szirmai,~J.:}
Isoptic curves of generalized conic sections in the hyperbolic plane.
\emph{Ukrainian Mathematical Journal}, {\bf 71/12} (2020), 1929-1944, doi: 10.1007/s11253-020-01756-3.
\bibitem{Ion}{Iona\c{s}cu, E. J.:} The "Circle" of Apollonius in Hyperbolic Geometry. \emph{Forum Geometricorum.}, Vol. {\bf 18}. (2018).
%
\bibitem{Cs-Sz25} Csima, G. -- Szirmai, J.: Translation-like Apollonius and triangular surfaces in non-constant curvature Thurston geometries. \textit{Results in Mathematics}, \textbf{80}, 190 (2025). https://doi.org/10.1007/s00025-025-02503-5
%
\bibitem{MSz}
{Moln{\'a}r,~E.~--~Szirmai,~J.:}
Symmetries in the 8 homogeneous 3-geometries.
\textit{Symmetry Cult. Sci.,}
{\bf 21/1-3}, 87-117 (2010).
%
\bibitem{Odehnal} Odehnal, B.: Equioptic curves of conic sections, 
\textit{J. Geom Graph.} \textbf{14} No.1 (2010), 29--43.
%
\bibitem{Sz} Szirmai, J.: Apollonius Surfaces, Circumscribed Spheres of Tetrahedra, Menelaus's and Ceva's Theorems in $\bS^2\!\times\!\bR$ and $\bH^2\!\times\!\bR$ Geometries. \textit{The Quarterly Journal of Mathematics}, \textbf{73.2} (2022), 477--494.
\bibitem{T} Thurston,~W.~P. (and Levy,~S. editor): Three-Dimensional Geometry and Topology.  Princeton University Press,  Princeton, New Jersey, vol. {\bf 1} (1997)
\end{thebibliography}
\end{document}